# An Essay on the Double Nature of the Probability

*Paolo Rocchi* *†;  *Leonida Gianfagna* *

**Abstract:** Classical statistics and Bayesian statistics refer to the frequentist and subjective theories of probability respectively. Von Mises and De Finetti, who authored those conceptualizations, provide interpretations of the probability that appear incompatible. This discrepancy raises ample debates and the foundations of the probability calculus emerge as a tricky, open issue so far.

Instead of developing philosophical discussion, this research resorts to analytical and mathematical methods. We present two theorems that sustain the validity of both the frequentist and the subjective views on the probability. Secondly we show how the double facets of the probability turn out to be consistent within the present logical frame.

## *1 – Introduction*

When Hilbert prepared the list of the most significant mathematical issues to tackle in the next future, he included the probability foundations in the group of 23 famous problems. Some years later Kolmorogov established the axioms of the probability calculus that the vast majority of authors have accepted so far. The Kolmogorovian axioms of non-negativity, normalization and finite additivity provide the rigorous base for calculations; they solve the issue in the terms posed by Hilbert but do not entirely unravel the fundamentals of the probability calculus. The probability is a measure and the physical nature of this measure has not been explained in a definitive manner [1]. The relation between probability and experience looks to be rather controversial.

For centuries the law of large numbers elucidated the meaning of the probability. The relative frequency in a large sample of trials offered the indisputable evidence for the probability calculated on the paper but after the First World War a significant turning point occurred. The calculus of the probability was used increasingly in business management. Many problems that arose in this environment dealt with making decisions under uncertainty and the method based on the frequency was not applicable in these cases as it is usually impossible to accumulate wide experience to assess the probability of an economical event. The subjective interpretation of the probability was introduced in this environment and became increasingly popular. In the first time theorists laid charge of arbitrariness against the subjective theory and vehemently rejected this frame [2], [3] and [4]. Later the Bayesian methods infiltrated professionals' environment and experimentalists learned to calculate the probability on the basis of acquired information [5].

However the pragmatic adoption of different statistical methods appears inconsistent due to the conceptual divergence emerging between the probability seen as a degree of belief and the probability as long-term relative frequency. Popper sums up this situation by the ensuing statement:

> "In modern physics (…) we still lack a satisfactory, consistent definition of probability; or what amounts to much the same, we still lack a satisfactory axiomatic system for the calculus of probability, [in consequence] physicists make much use of probabilities without being able to say, consistently, what they mean by 'probability'" [6]

This paper is an attempt to show why and when the frequentist and subjective approaches can be considered compatible in point of logic. We mean to prove that the dualist view on probability turns out to be consistent using two theorems.

---

* IBM, via Shangai 53, Roma, Italy
† Luiss University, via Alberoni 7, Roma, Italy

## 2. Two Theorems

We calculate some properties of the probability in relation to the relative frequency under the following assumptions:

**1)** The relative frequency is the observed number $q$ of successful events $A$ for a finite sample of $n$ trials

$$F(A_n) = q/n \tag{1}$$

We assume $F(A_n)$ as the counterpart of $P(A)$ in the real world.

**2)** We exclude that $A$ is sure or impossible, and $X_1, X_2, X_3,.. X_m$ are the possible results of $A$ with finite *expected value* $\mu$

$$\mu = \frac{X_1 + X_2 + X_3 + ... + X_m}{m} \tag{2}$$

In principle there are two extreme possibilities:

**3)** There are $n$ trials with

$$n > 1 \tag{3}$$

The occurrences of $A$ make a sequence of pair-wise independent, identically distributed random variables $X_1, X_2, X_3,.. X_n$ with finite sample mean $SM_n$

$$SM_n = \frac{1}{n}\sum_{j=1}^{n} X_j . \tag{4}$$

**4)** There is only one trial namely

$$n = 1 \tag{5}$$

### 2.1 Theorem of large numbers (strong)

Suppose hypothesis **3)** true, then we have that $\bar{X}_n$ converges almost surely to $\mu$

$$P\left(\lim_{n\to\infty} SM_n = \mu\right) = 1 \tag{6}$$

*Proof*: There are several methods to prove this theorem. We suggest the recent proof by Terence Tao who reaches the result with rigor under the mathematical framework of measure theory [7]

*Remarks*: This theorem expresses the idea that as the number of trials of a random process increases, the percentage difference between the actual and the expected values goes to zero. It follows from this theorem that the relative frequency of success in a series of Bernoulli trials will converge to the theoretical probability

$$F(A_n) \to P(A) \quad \text{when } n \to \infty \tag{7}$$

The relative frequency $F(A_n)$ after $n$ trials will almost surely converge to the probability $P(A)$ as $n$ approaches infinity.

### 2.2 Theorem of a single number

Let **4)** true, the relative frequency of the event is not equal to the probability of $A$

$$F(A_1) \neq P(A) \tag{8}$$

*Proof:* The random event $A$ either occur or does not occur in a sole trial, the following values are allowed for the frequency

$$F(A_1) = 1/1 = 1 \tag{9}$$
$$F(A_1) = 0/1 = 0 \tag{10}$$

We have excluded that $A$ is sure or impossible, thus (8) and (9) mismatch with $P(A)$ that is decimal.

## 3 Analytical Approach

Speaking in general, the abstract determination of the parameter $\mathfrak{M}(\mathfrak{A})$ that refers to the argument $\mathfrak{A}$ is neutral, namely pure mathematicians have no concern for the material qualities of $\mathfrak{M}$ and of $\mathfrak{A}$; instead those tangible properties are essential for sperimentalists who particularize the significance of $\mathfrak{M}(\mathfrak{A})$ determined by $\mathfrak{A}$. As an example take the real function $y = f(x)$ that has immutable mathematical properties in the environments A, B and C beyond any doubt. Now suppose that the variables $x_A$, $x_B$ and $x_C$ are heterogeneous in physical terms, than the qualities of $y_A$, $y_B$ and $y_C$ cannot be compared. Even if one adopts identical algorithms to obtain $y_A$, $y_B$ and $y_C$, the concrete features of $y_A$, $y_B$ and $y_C$ do not show any analogy.

Kolmogorov chooses the subset $A$ belonging to the set space $\Omega$ as the flexible and appropriate model for the argument of the probability $P(A)$. Pure mathematicians usually assume that $A$ is the known term of a problem and do not call into question the subset $A$. This behavior works perfectly until one argues on the theoretical plane but things go differently in applications. We have discussed a set of misleading and paradoxes depending on the poor analysis of the probability argument [8], [9]. In this paper we dissect the physical characteristics of the probability through its argument.

Hypotheses (3) and (5) deal with the occurrences of a random event; and the subsets $A_n$ and $A_1$ prove to have two special states of being present in the world. No doubt $P(A_1)$ and $P(A_n)$ are equal as numbers but $A_1$ and $A_n$ turn out to be distinct arguments in Section 2. The theorems just demonstrated show that the probability has different applied properties depending on the argument, and we rewrite (7) in more precise terms:

$$F(A_n) \to P(A_n) \qquad \text{when } n \to \infty \tag{11}$$

And result (8) becomes:

$$F(A_1) \neq P(A_1) \tag{12}$$

The distinction between $A_n$ and $A_1$ is not new. Ramsey, De Finetti and other subjectivists calculate the probability of the single occurrence $A_1$ in contrast with Von Mises who coined the term '*collective*' to emphasize the large number of occurrences pertaining to $A_n$. Various dualists argue upon these extreme possibilities such as Popper who develops the *long-run propensity theory* associated with repeatable conditions, and the *single-case propensity theory*.

It is necessary to note that philosophers, who argue on the probabilities of single and multiple events, are inclined to integrate the probability together with its argument. Generally speaking philosophers search for a comprehensive solution, they melt the probability with its argument and obtain a synthetic theoretical notion. Instead the present inquiry is analytical. We keep apart $A_n$ from $A_1$ and later deduce the features of $P$ from each argument. To exemplify statement (11) specifies the property of $A_n$ instead one cannot tell the same for the argument $A_1$. Theorem 2.1 does not provide conclusions for the probability of a single trial, but solely for $P(A_n)$ in that $P(A_n)$ depends on $A_n$.

### 3.1 When the Number is Large

Result (11) details that the relative frequency obtained after repeated trials in a Bernoulli series conforms to the probability. Because of this substantial experiment one can conclude that $P(A_n)$ is a physical measurable quantity

$$P(A_n) \quad \text{is a real quantity in the world.} \tag{13}$$

The weak point of theorem 2.1 is the limit conditions for the validation of $P(A_n)$. Solely infinite trials can corroborate the quantity $P(A_n)$ and nobody can bring forth this ideal experiment. If one does not reason in abstract, he can conclude that the ideal test does not place insurmountable obstacles. Experimentalists are aware that the higher is $n$ and the better is the achieved validation of $P(A_n)$. The larger is the number of trials, the more the probability finds solid confirmation namely one has the rigorous method to validate the probability in the world out there.

Theorem of large numbers has long-term prehistory. Since immemorial time gamblers and laymen had keen interest in the corroboration of the values obtained through previsions, and they sensed that if an experiment is repeated a large number of times then the relative frequency with which an event occurs tends to the theoretical possibilities of that event.

Ample literature deals with the large number theorem and further comments seem unnecessary to us. Instead the single number theorem requires more accurate annotations.

## 3.2 When the Number is the Unit

A premise turns out to be necessary to the discussion of inequality (12).

Any hypothesis must be supported by appropriate validation in sciences. Theories and assertions are to be confirmed by evidences and not just by mental speculation [10]. Scientists search, explore, and gather the evidences that can sustain or deny whatsoever opinion. Popper highlights how a sole experiment can disprove a theoretical frame. Many candidate interpretations have been proposed by theoretical physicists during the centuries that have not been confirmed experimentally.

Test is supreme and a quantity that cannot be verified through trials is not extant. For examples some geologists suppose that earthquakes are able to emit premonitory signals. Forewarning signs could save several lives, but accurate tests do not bring indisputable evidences thus the scientific community concludes that premonitory signals do not exist.

Eqn. (12) holds that the frequency $F(A_1)$ systematically mismatches with $P(A_1)$ and this entails that one cannot validate the probability of $A_1$ in the physical reality. It is not a question of operational obstacles, errors of instruments, inaccuracy or other forms of limitation. Statement (12) establishes that in general one cannot corroborate the probability $P(A_1)$ by means of tests thus theorem 2.2 leads to the following conclusion:

$$P(A_1) \quad \text{is a quantity that does not exist in the world.} \tag{14}$$

The probability of a single event is a number written on the paper which none can control in the real world.

This conclusion turns out to be dramatic as statistics play crucial role in sciences and modern economies. Statistics offers the necessary tools to investigate and understand complex phenomena. Moreover people have strong concern in unique events whose outcomes sometimes are essential for individuals to survive. In order to not throw away $P(A_1)$ authors established to bypass (14) in the following way.

The single-case probability has symbolic value, it is a sign and as such carries out communication task. $P(A_1)$ is a piece of information that can influence the human soul and can sustain personal credence. Thus the probability $P(A_1)$ may be used to express the degree of belief of an individual about the outcome of an uncertain even. Theorists do not discard the probability $P(A_1)$ and ground the subjective theory of probability upon this original function assigned to $P(A_1)$.

In principle subjective probability could vary from person to person, and subjectivists answer back the accusation of arbitrariness. They do not deem 'the degree of belief" as something measured by strength of feeling, but in terms of rational and coherent behavior. The "Dutch book" criterion establishes a guide for the agent's degrees of belief which have to satisfy synchronic coherence conditions. An agent evaluates the plausibility of a belief according to past experiences and to knowledge coming from various sources, and this criterion opens the way to the Bayesian methods [11].

Theorem 2.2 has significant antecedents.

For centuries the vast majority of philosophers shared the deterministic logic and rejected any form of probabilistic issues. This orientation emerged in the Western countries since the classical era. In fact thinkers founded out that the possible results of a single random event yield unrealistic conclusions. As an example, suppose to flip a coin, one forecasts two balanced possibilities, instead a single test exhibits head or tail. Indeterminist logic suggests two outcomes equally likely, conversely experience shows only one outcome. Philosophers deemed a logic topic not entirely determined as unreliable and in consequence of this severe judgment the apodictic reasoning dominated the scientific realm. This systematic orientation kept scientists from paying adequate attention to probabilistic problems, and statistical methods infiltrated scientific environments solely by the end of the nineteenth century [12].

One can apply the frequentist theory provided the set of trials is very large instead the Bayesian statistics works even for a sole trial. The subjective theory does not encounter the obstacles typical of the frequentist frame. De Finetti and others recognize the theorem of large number as a property of the single event that repeats several times and conclude that the subjective theory is universal [13].

In consequence of (13) and (14) we can answer in the following terms.

Theorem 2.1 holds that $P(A_n)$ is a real parameter which one can validate in the world; conversely theorem 2.2 entails that $P(A_1)$ translates a personal belief into a number. The indefinite series $P_1(A_1), P_2(A_1), P_3(A_1), P_4(A_1)$... mirrors personal credence upon the repeated event $A_1$ and cannot be compared to $P(A_n)$ which is an authentic measure controllable in the physical reality. Hence it is worth the ponderous verification of $P(A_n)$ since this is a physical measure; instead $P(A_1)$ is a symbol.

It is true that $P(A_1)$ is flexible to countless single situations, but the theorem of a single number proves that $P(A_1)$ does not have any counterpart in the physical reality.

## 4. The Nature of the Probability

Least we can answer from theorems 2.1 and 2.2 whether the probability is a measure with physical meaning, or is not physical or even has a double nature.

Theorem 2.1 demonstrates that $P(A_n)$ have different properties respect to $P(A_1)$. Statements (13) and (14) hold that the frequentist nature and the subjective nature of the probability relay on two separated situations. Therefore (13) and (14) prove to be compatible in point of logic since $P(A_1)$ and $P(A_n)$ refer to distinct physical events. It may be said that $P(A)$ is a real, controlled measure when:

$$A = A_n \qquad (15)$$

Instead $P(A)$ has a personal value when:

$$A = A_1 \qquad (16)$$

It is reasonable to conclude that the probability is double in nature and professional practice sustains this conclusion as the classical statistics and the Bayesian statistics investigate different kinds of situations. They propose methods that refer to separated environments. To exemplify physicists who validate a general law adopt the Fisher statistics; instead when a physicist has to make decision on a single experiment usually he adopts the Bayesian methods.

## 5. Conclusion

The interpretation of the probability is considered a tricky issue since decades, and normally theorists tackle this broad argument from the philosophical stance. The philosophical approach yielded several intriguing contributions but did not provide the definitive solution so far. We fear philosophy is unable to bring forth the overwhelming proof upon the nature of the probability because of its all-embracing method of reasoning. Instead the present paper suggests the analytical approach which keeps apart the concept of probability from its argument and examines the influence of $A_1$ and $A_n$ over $P$.

Nowadays some statisticians are firm adherents of one or other of the statistical schools, and endorse their adhesion to a school when they begin a project; in a second stage each expert uses the statistical method which relies on his opinion. However this faithful support for a school of thought may be considered an arbitrary act as relying on personal will.
Another accepted view is that each probabilistic theory has strengths and weaknesses and that one or the other may be preferred. This generic behavior appears disputable in many respects since one overlooks the conclusions expressed by Von Mises which clash against De Finetti's conceptualization.
The present frame implies that the use of a probability theory is not a matter of opinion and suggests a precise procedure to follow in the professional practice on the basis of incontrovertible facts. Firstly an expert should see whether he is dealing with only one occurrence or with a high number of occurrences. Secondly an expert should follow the appropriate method of calculus depending on the extension (15) or (16) of the subset $A$.

The present research is still in progress. This paper basically focuses on classical probability and we should extend this study to quantum probability in the next future.

## References

[1] Gillies D. – *Philosophical Theories of Probability* – Routledge (2000)
[2] Mayo D. G. – *Error and the Growth of Experimental Knowledge* – University of Chicago Press (1996).
[3] Forster M. R. – How Do Simple Rules "Fit to Reality" in a Complex World? – *Minds and Machines*, 9, pp. 543-564 (1999).
[4] Denker M., Woyczyński W.A., Ycar B. – Introductory Statistics and Random Phenomena: Uncertainty, Complexity and Chaotic Behavior in Engineering and Science – Springer-Verlag Gmbh (1998)
[5] D'Agostini G. – Role and Meaning of Subjective Probability; Some Comments on Common Misconceptions – *XX Intl. Workshop on Bayesian Inference and Maximum Entropy Methods in Science and Engineering,* AIP Conference Proceedings, vol. 568, 23 (2001).
[6] Popper K. R. - *The Logic of Scientific Discovery* – Routledge (2002).
[7] Tao T. – *The Strong Law of Large Numbers* – http://terrytao.wordpress.com/2008/06/18/the-strong-law-of-large-numbers/ (accessed 15 January 2012).
[8] Rocchi P., Gianfagna L. – Probabilistic Events and Physical Reality: a Complete Algebra of Probability - *Physics Essays*, 15(3), pp. 331-118 (2002). http://arxiv.org/ftp/math/papers/0309/0309069.pdf
[9] Rocchi P. – *The Structural Theory of Probability* – Kluwer/Plenum (2003).

[10] Gower B. – *Scientific Method: A Historical and Philosophical Introduction* – Taylor & Francis (2007).
[11] Savage L. J. – The Foundations of Statistics – Courier Dover Publications (1972).
[12] Earman J. – Aspects of Determinism in Modern Physics – in *Philosophy of Physics, Part B,* Butterfield J., Earman J. (eds.) North Holland, pp. 1369-1434 (2007).
[13] De Finetti B. – *Theory of Probability: A Critical Introductory Treatment* – John Wiley & Sons (1975).